\documentclass[runningheads]{llncs}
\usepackage{comment}
\usepackage{graphicx}
\usepackage{hyperref}
\usepackage{amsmath,amssymb,amsfonts}%
\usepackage{mathrsfs}%
\usepackage{xcolor}%
\usepackage{textcomp}%
\usepackage{algorithm}%
\usepackage{algorithmicx}%
\usepackage{algpseudocode}%
\usepackage{listings}%
\usepackage{mathtools}
\usepackage{subcaption}
\begin{document}
\title{Krylov Solvers for Interior Point Methods with Applications in Radiation Therapy and Support Vector Machines}
\titlerunning{Krylov Solvers for Interior Point Methods}

\author{Felix Liu \inst{1,2}\orcidID{0000-0001-6865-9379} \and
Albin Fredriksson \inst{2}\and
Stefano Markidis \inst{1}}

\institute{KTH Royal Institute of Technology, Stockholm, Sweden \\
\email{felixliu@kth.se} \and 
RaySearch Laboratories, Stockholm, Sweden}
\authorrunning{F. Liu et al.}
%
%
\maketitle              
\begin{abstract}
Interior point methods are widely used for different types of mathematical optimization problems. Many implementations of interior point methods in use today rely on direct linear solvers to solve systems of equations in each iteration. The need to solve ever larger optimization problems more efficiently and the rise of hardware accelerators for general purpose computing has led to a large interest in using iterative linear solvers instead, with the major issue being inevitable ill-conditioning of the linear systems arising as the optimization progresses. We investigate the use of Krylov solvers for interior point methods in solving optimization problems from radiation therapy and support vector machines. We implement a prototype interior point method using a so called doubly augmented formulation of the Karush-Kuhn-Tucker linear system of equations, originally proposed by Forsgren and Gill, and evaluate its performance on real optimization problems from radiation therapy and support vector machines. Crucially, our implementation uses a preconditioned conjugate gradient method with Jacobi preconditioning internally. Our measurements of the conditioning of the linear systems indicate that the Jacobi preconditioner improves the conditioning of the systems to a degree that they can be solved iteratively, but there is room for further improvement in that regard. Furthermore, profiling of our prototype code shows that it is suitable for GPU acceleration, which may further improve its performance in practice. Overall, our results indicate that our method can find solutions of acceptable accuracy in reasonable time, even with a simple Jacobi preconditioner.

\keywords{Interior point method \and Krylov solver \and Radiation therapy \and Support Vector Machines}
\end{abstract}
\section{Introduction}
Mathematical optimization is used in many areas of science and industry, with applications in fields like precision medicine, operations research and many others. In this paper, we focus on the solution of quadratic programs (QP), which are optimization problems with a quadratic objective function and linear constraints, using an interior point method (IPM). These arise in many applications naturally, for instance in training support vector machine classifiers, but can also be used as part of a sequential quadratic programming solver \cite{boggs1995sequential} to solve more general nonlinear optimization problems.


Computationally, IPMs for optimization rely on Newton's method to find search directions. This involves the solution of a large, often sparse and structured linear system of equations, which we will refer to henceforth as the Karush-Kuhn-Tucker (KKT) system, at each iteration. Traditionally, this system is often solved using direct linear solvers, such as $LDL^T$-factorization for indefinite matrices or Cholesky factorization for positive definite methods, but a topic of interest for much research in the field is the use of iterative linear solvers \cite{saad2003iterative} instead. Indeed, the move to iterative linear algebra for interior point methods has been identified by some authors as a key step in enabling interior point methods to handle very large optimization problems \cite{gondzio2012interior}.

Another advantage of iterative algorithms for solving linear systems is their suitability for modern computing hardware, such as GPUs. With their rising dominance in High-Performance Computing (HPC), extracting the maximum performance from modern computing hardware all but requires the use of some type of accelerator. Direct linear solvers can suffer performance wise on these types of hardware for a variety of reasons, e.g., unstructured memory accesses, which has been seen in previous studies \cite{swirydowicz2022linear} in the context of interior point methods. The major challenge of using iterative solvers lies in a structured form of ill-conditioning that inevitably exists in the linear systems, which can be severely detrimental to the convergence of the linear solver. Still, given potential performance gains from the use of significantly more powerful computing resources, we believe the trade-off between numerical stability and parallel performance is worthy of further investigation.

One type of problems we consider arise from treatment planning for radiation therapy, which loosely speaking is the process of optimizing treatment plans (treatment machine parameters) for each individual patient case to deliver as accurate a dose as possible to the tumor volume. This inverse problem is often solved by formulating it as a constrained optimization problem, for which finding a solution can be both computationally expensive and present an important bottleneck in the clinical workflow. Computational speed and efficiency is thus of crucial importance. To demonstrate the applicability of our proposed method to other problems as well, we also consider problems from the training of support vector machine classifiers \cite{noble2006support}, an important method from classical machine learning. All in all, we are interested in studying and evaluating the potential of using IPMs with iterative linear solvers as an avenue to enable us to utilize accelerators and powerful computing resources for solving these problems more efficiently.

In this paper, we propose a complete IPM solver prototype based on the work by Forsgren and Gill in \cite{forsgren2007iterative}, where a special formulation of the KKT systems from interior point methods is considered, which guarantees that the KKT system is positive-definite throughout the optimization for convex problems. Our contribution a complete IPM solver prototype using the doubly augmented formulation and iterative linear solvers, which is capable of solving real world optimization problems. We demonstrate its effectiveness for the solution of quadratic optimization problems arising from real world applications in both radiation therapy treatment planning as well as support vector machines.
\section{Background}
\subsection{Interior Point Methods}
We are concerned with the solution of convex, continuous quadratic programs using interior point methods. The following section introduces the relevant aspects for our proposed method. Readers interested in a more thorough overview of interior point methods for optimization are referred to e.g. \cite{wright1992interior,forsgren2002interior}. In general, we will be dealing with a problem on the form:
\begin{equation}
\begin{aligned}
    \text{min.}& \quad \frac{1}{2} x^T H x + p^T x \\
    \text{s.t.}& \quad l \leq Ax \leq u \\
               & \quad C x = b
\end{aligned}
\end{equation}
where $H$ is the (positive definite) Hessian of the objective function, $p \in \mathbb{R}^n$ are the linear coefficients of the objective, $A, C$ are the Jacobians of the (linear) inequality- and equality constraints, respectively. In general, we allow components of the constraints to be unbounded, but we do not account for this explicitly to simplify the exposition. Inequality constraints are often treated by the introduction of auxiliary \emph{slack variables}. We convert the problem to having only lower bounds, introduce slack variables $s_l, s_u$ (for lower and upper bounds respectively), use a log-barrier term for the inequality constraints and finally, we use a penalty barrier method \cite{fiacco1990nonlinear,forsgren2002interior} to handle equality constraints, giving:
\begin{equation}
\begin{aligned}
    \text{min.} \quad &\frac{1}{2} x^T H x + p^T x - \\ 
                      &- \mu \sum \log((s_l)_i) - \mu \sum \log((s_u)_i) + \\
                      &+ \frac{1}{2\mu} ||Cx - b||^2\\
    \text{s.t.} \quad &Ax - s_l - l = 0 \\
                      -&Ax - s_u + u = 0\\
\end{aligned}
\label{eq:qp_barrier}
\end{equation}
Here $\mu$ is the so called barrier parameter. Intuitively, the barrier terms diverge towards $+\infty$ as the boundary of the feasible set is approached, thus encouraging feasibility throughout the iterations. As is common in primal-dual interior point methods, we work with the \emph{perturbed} optimality conditions. These state that an optimal solution to the equality constrained subproblem in \eqref{eq:qp_barrier} satisfies the following:
\begin{equation}
\label{eq:perturbed_KKT}
\begin{aligned}
    &r_H \coloneqq Hx + p - A^T \lambda_l + A^T \lambda_u - C^T \lambda_e = 0 \\
    &r_{A_l} \coloneqq Ax - s_l - l = 0 \\
    &r_{A_u} \coloneqq -Ax + s_u + u = 0 \\
    &r_e     \coloneqq Cx - b + \mu \lambda_e = 0\\
    &r_{c_1} \coloneqq (\lambda_l)_i (s_l)_i - \mu = 0, \quad i = 1,...,m_l \\
    &r_{c_2} \coloneqq (\lambda_l)_i (s_u)_i - \mu = 0, \quad i = 1,...,m_l.\\
\end{aligned}
\end{equation}
These conditions are very similar to the first order Karush-Kuhn-Tucker conditions for optimality \cite[Ch. 12.3]{nocedal2006numerical}, except that the final two conditions are perturbed by $\mu$, and the inclusion of the penalty barrier method for equality constraints. Primal-dual interior point methods generally seek points satisfying the perturbed optimality conditions above using, e.g., Newton's method, while successively decreasing the barrier parameter $\mu \rightarrow 0$. As $\mu$ approaches $0$, we expect our solution to approach a point satisfying the KKT conditions for optimality.
\subsection{Optimization Problems in Radiation Therapy and SVMs} \label{sec:radiotherapy_opt}
The first type of problem we consider are from radiation therapy, and are all exported from the treatment planning system RayStation, developed by the Stockholm-based company RaySearch Laboratories. The problems all arise from \emph{treatment planning} for radiation therapy, where an optimization problem is solved to determine a \emph{treatment plan} for each individual patient.
A view of treatment planning from an optimization perspective can be found in e.g. \cite{ehrgott2010mathematical,engberg2018automated}.

The optimization problems from RayStation are QP-subproblems in a Sequential Quadratic Programming solver used for nonlinear optimization, and have the form:
\begin{equation}
\begin{aligned}
    \text{min.} \quad &\frac{1}{2} p^T \nabla_{xx} \mathcal{L}(x, \lambda) p + (\nabla f)^T(x) p \\
    \text{s.t.} \quad &\nabla g(x)^T p + g(x) \geq 0.
\end{aligned}
\end{equation}
Here, $\mathcal{L}(x, \lambda) = f(x) - g(x) \lambda^T$ is the Lagrangian, and $\lambda$ are the Lagrange multipliers. In practice, the Hessian of the Lagrangian $\nabla_{xx} \mathcal{L}(x, \lambda)$ can be expensive to form, and it is common to use a quasi-Newton type approximation of the Hessian instead. The SQP solver uses Broyden-Fletcher-Goldfarb-Shanno (BFGS) updates \cite{broyden1970convergence} to estimate the Hessian of the non-linear problem, which means that the Hessian for each of our QP-subproblems can be written on matrix form as:
\begin{equation}
    H = H_0 + UWU^T,
\end{equation}
where $H_0$ is the initial approximation to the Hessian, the dense $n \times k$ matrix $U$ consists of the update vectors on the columns, and $W$ is a diagonal matrix with the scalar update weights on the diagonal. With a suitable line-search method, the updates to the Hessian can be ensured to preserve positive definiteness, thus making the QP-subproblem we need to solve at each SQP iteration convex.

The second type of problem we consider are QP dual problems from support vector machine training for classification. These problems are of the form:
\begin{equation}
\begin{aligned}
    \text{min.} \quad &\frac{1}{2} \alpha^T H \alpha - \alpha^T e \\
    \text{s.t.} \quad &\alpha^T y = 0 \\
                      &0 \leq \alpha \leq c,
\end{aligned}
\end{equation}
where the Hessian $H$ is of the form $H = y y^T Q$, and the entries of $Q$ are $K(x_i, x_j)$, for some "kernel" $K$ used to map the data into a (more) separable space. For our experiments, we use a radial basis function kernel
\begin{equation*}
K(x_i, x_j) = \exp(-||x_i - x_j||^2 / 2\sigma).
\end{equation*}
\section{Prototype Interior Point Method}
We now describe the key components of our prototype interior point method implementation used in this paper.
\subsection{KKT System Formulation}
Many optimization problems from real applications include bounds on the (primal) variables themselves. Such bounds can be included by the introduction of appropriate rows in the $A$ matrix. For an efficient formulation, we will consider the bounds on variables separately from general linear constraints.
Separating the handling of the variable bounds (i.e. lower and upper bounds on the values of the variables $x$) and then using Newton's method to solve the perturbed optimality conditions \eqref{eq:perturbed_KKT} gives a linear system of the form:
\begin{equation}
    \begin{pmatrix}
        H & -A^T & A^T & -I & I & C^T &    &    &    & \\
        C  &     &     &    &   & M &    &    &    & \\
        A  &     &     &    &   &     & -I &    &    & \\
        -A &     &     &    &   &     &    & -I &    & \\
        I  &     &     &    &   &     &    &    & -I & \\
        -I &     &     &    &   &     &    &    &    & -I \\
           & S_{l_A} & &    &   &     &\Lambda_{l_A} & & & \\
           & & S_{u_A} &    &   &     &  &\Lambda_{u_A}  & & \\
           & & & S_{l_x} &  &   &     &  & \Lambda_{l_x} & \\
           & & & & S_{u_x}  &   &     &  &  &  \Lambda_{u_x}
    \end{pmatrix}
    \begin{pmatrix}
        \Delta x \\
        \Delta \lambda_{l_A} \\
        \Delta \lambda_{u_A} \\
        \Delta \lambda_{l_x} \\
        \Delta \lambda_{u_x} \\
        \Delta \lambda_e \\
        \Delta s_{l_A} \\
        \Delta s_{u_A} \\
        \Delta s_{l_x} \\
        \Delta s_{u_x} \\
    \end{pmatrix}
    =
    \begin{pmatrix}
        -r_H \\
        -r_e \\
        -r_{l_A} \\
        -r_{u_A} \\
        -r_{l_x} \\
        -r_{u_x} \\
        -r_{c_1} \\
        -r_{c_2} \\
        -r_{c_3} \\
        -r_{c_4} \\
    \end{pmatrix},
\label{eq:kkt_full}
\end{equation}
where $\lambda_e, \lambda_{l_A}, \lambda_{u_A}, \lambda_{l_x}, \lambda_{u_x}$ are the Lagrange multipliers for the equality constraints, lower and upper bounds on the (general) linear constraints, and lower and upper bounds for the variable bounds respectively. The slack variables $s$ are subscripted in the same way. The residuals on the RHS are similar to the ones shown in \eqref{eq:perturbed_KKT}. $\Lambda, S, M$ denote diagonal matrices with the corresponding Lagrange multipliers, slack variables or barrier parameter $\mu$ on the diagonal respectively, and $e$ is an appropriately sized column-vector with a value of 1 in all coefficients.

The above system can be reduced in size by block-row elimination. Multiplying the sixth and seventh block row by $\Lambda_{l_A}^{-1}$ and the fifth block row by $\Lambda_{u_A}^{-1}$ and adding them to the second and third rows, followed by multiplying the sixth and seventh block rows by $S_{l_x}^{-1}$ and $-S_{u_x}^{-1}$ respectively and adding to the top row, as well as multiplying the fourth and fifth block rows by $S_{l_x}^{-1}\Lambda_{l_x}$ and $-S_{u_x}^{-1}\Lambda_{u_x}$ and adding to the top row. The final reduced linear system of equations (with the same row operations on the RHS) can be written as:
\begin{equation}
    \begin{pmatrix}
    Q & -B^T \\
    B & D
    \end{pmatrix}
    \begin{pmatrix}
        \Delta x \\
        \Delta \lambda_{A}
    \end{pmatrix}
    =
    \begin{pmatrix}
        r_1 \\
        r_2
    \end{pmatrix},
    \label{eq:2x2_augmented}
\end{equation}
where:
\begin{equation*}
\begin{aligned}
    &Q = H + S_{l_x}^{-1} \Lambda_{l_x} + S_{u_x}^{-1} \Lambda_{u_x}, \quad
    B = \begin{pmatrix}
        C \\
        A \\
        -A
    \end{pmatrix}, \quad
    D = \begin{pmatrix}
        M &  &  \\
          & \Lambda_{l_A}^{-1} S_{l_A} &  \\
          & & \Lambda_{u_A}^{-1} S_{u_A}
    \end{pmatrix}. \\
    &\Delta \lambda_A =
    \begin{pmatrix}
        \Delta \lambda_e     \\
        \Delta \lambda_{l_A} \\
        \Delta \lambda_{u_A}
    \end{pmatrix}, \quad
    r_1 = -r_H - S_{l_x}^{-1} r_{c_3} + S_{u_x}^{-1} r_{c_4} - S_{l_x}^{-1} \Lambda_{l_x} r_{l_x} + S_{u_x}^{-1} \Lambda_{u_x} r_{u_x}\\
    &r_2 =
    \begin{pmatrix}
        -r_e \\
        -r_{A_l} - \Lambda_{l_A}^{-1} r_{c_1} \\
        -r_{A_u} - \Lambda_{u_A}^{-1} r_{c_2}
    \end{pmatrix}
\end{aligned}
\end{equation*}
At this point the purpose of handling the variable bounds separately becomes clear, since they now only contribute a diagonal term $S_{l_x}^{-1} \Lambda_{l_x} + S_{u_x}^{-1} \Lambda_{u_x}$ in the Hessian block. A challenge with the system \eqref{eq:2x2_augmented} is that it becomes inevitably ill-conditioned as the optimization approaches a solution. Intuitively this can be seen by noting that as $\mu \rightarrow 0$, some elements of the diagonal block $D$ become very small and some become unbounded, since for active constraints the slack variables tend to zero, while for inactive constraints the Lagrange multipliers do. For more details on this ill-conditioning see, e.g., \cite{wright1998ill,forsgren1996stability}.
The system \eqref{eq:2x2_augmented} is unsymmetric, and it is common to consider many equivalent but symmetric systems instead. Our implementation uses the \emph{doubly augmented} formulation, proposed in \cite{forsgren2007iterative}. This formulation can be derived through block-row operations on the system in \eqref{eq:2x2_augmented}, by multiplying the second block row by $2B^T D^{-1}$ and adding it to the first block row:
\begin{equation}
    \begin{pmatrix}
    Q + 2B^T D^{-1} B & B^T \\
    B & D
    \end{pmatrix}
    \begin{pmatrix}
        \Delta x \\
        \Delta \lambda_{A} \\
    \end{pmatrix}
    =
    \begin{pmatrix}
        r_1 + 2B^T D^{-1} r_2\\
        r_2
    \end{pmatrix}.
    \label{eq:doubly_augmented}
\end{equation}
The major advantage is that the matrix is symmetric and positive definite for convex problems \cite{forsgren2007iterative}. This enables us to use a preconditioned conjugate gradient method to solve the system efficiently. To precondition the system we use Jacobi (diagonal scaling) preconditioning, which is motivated by the ill-conditioning arising primarily from the diagonal $D$ block in the matrix.
\subsection{IPM Implementation}
We implement a prototype interior point method to assess the performance and accuracy of the method when using iterative linear algebra. The following is a brief description of the design of our implementation.
Our implementation is a primal-dual interior point method based on the KKT system formulation in \eqref{eq:doubly_augmented}. We use a ratio test to determine the maximum step length to take each iteration to maintain the positivity of the slack variables and Lagrange multipliers:
\begin{equation}
\begin{aligned}
    &\alpha_x = \min \left\{1.0, \gamma \left(\min\left\{ -\frac{s_i}{\Delta s_i} : \Delta s_i < 0 \right\}\right)\right\} \\
    &\alpha_{\lambda} = \min \left\{1.0, \gamma \left(\min\left\{ -\frac{\lambda_i}{\Delta \lambda_i} : \Delta \lambda_i < 0 \right\}\right)\right\}.
\end{aligned}
\label{eq:line_search}
\end{equation}
The scalar $\gamma < 1$ is to ensure strict positivity of the slacks and Lagrange multipliers throughout the optimization, and we use a value of $\gamma = 0.99$ in our implementation. The step lengths from the line search are used to scale the search direction.

Finally, we decrease the barrier parameter $\mu$ based value of the residuals (shown in the right-hand side of \eqref{eq:kkt_full}). Namely, when the 2-norm of the residuals is smaller than the current value of $\mu$, we divide $\mu$ by 10 and continue the optimization. The optimizer terminates when the barrier parameter decreases below a tolerance threshold, which by default is set to $10^{-6}$.

We summarize the main components of our implementation in Algorithm~\ref{alg:IPM_prototype}.
\begin{algorithm}[h!]
\caption{Interior Point Method}\label{alg:IPM_prototype}
\begin{algorithmic}[1]
    \For{$i \gets 1$ to $N$}
        \State Find search direction by solving \eqref{eq:2x2_augmented} using PCG
        \State Line search for $\alpha_x, \alpha_{\lambda}$ from \eqref{eq:line_search}
        \State $x \gets x + \alpha_x \Delta x$
        \State $\lambda \gets \lambda +  \alpha_{\lambda} \Delta \lambda$
        \State $s \gets s + \alpha_x \Delta s$
        \State Update diagonal $D$ in KKT system
        \State Compute residuals (RHS of \eqref{eq:kkt_full})
        \If{$||r|| < \mu$} \Comment{$r$ is the RHS of \eqref{eq:kkt_full}}
            \If{$\mu < \mu_{tol}$}
                \State Return solution
            \EndIf
            \State $\mu \gets \mu / 10$
        \EndIf
    \EndFor
\end{algorithmic}
\end{algorithm}
The use of a Krylov solver for our implementation provides some practical benefits in how the computation is structured, as it allows us to work in a matrix-free manner. Concretely, we do not explicitly form the $2B^TD^{-1}B$ term in \eqref{eq:2x2_augmented} (nor the entirety of the matrix), but always work with the different components separately. This is also especially advantageous for the quasi-Newton structure of the Hessian. as discussed in Section~\ref{sec:radiotherapy_opt}. Recall that the BFGS-Hessian can be written in matrix form as $H = H_0 + UWU^T$, which is a dense $n \times n$ matrix, where $n$ is the number of variables in the QP. Similarly to before, this matrix is also not explicitly formed in our solver, saving significant computational effort when the number of variables is large.

We have implemented the method described in C++ (using BLAS for many computational kernels), which is also the implementation we use for the experiments conducted in this work.
\section{Experimental setup}
\begin{table}[]
    \centering
    \caption{Problem dimensions for the considered QPs}
    \begin{tabular}{|c|c|c|c|}
        \hline
        Problem & Vars. & Lin. cons. & Bound cons.\\
        \hline
        \rule{0pt}{3ex}    
        Proton H\&N & 55770 & 0 & 90657 \\
        Proton Liver & 90657 & 15 & 90657 \\
        Photon H\&N & 13425 & 42273 & 13425 \\
        SVM a1a & 1605 & 1 & 3210  \\
        \hline
    \end{tabular}
    \label{tab:problems}
\end{table}
\noindent
The radiation therapy optimization problems evaluated in this work are exported directly from the RayStation optimizer. We export the QP subproblems to files directly from RayStation which permits us to use them for our experiments, without relying on RayStation itself. In particular, we consider three problem cases, two cases treated using proton therapy, one for the head-and-neck region and a liver case and one case treated using photons. The dimensions of the corresponding optimization problems are shown in Table~\ref{tab:problems}. In all considered problems, the problems are exported from the later stages of the SQP iterations, which are typically the most challenging.

For the SVM training problem, we use the \texttt{a1a} problem available from the LIBSVM dataset \cite{chang2011libsvm}. We pre-compute the (dense) Hessian H using the radial basis function kernel as described in Section~\ref{sec:radiotherapy_opt}.

The performance measurements were all carried out on a local workstation equipped with an AMD Ryzen 7900x CPU with 64 GB of DDR5 DRAM. The BLAS library used was OpenBLAS 0.3.21 with OpenMP threading. The measurements of the condition number were run on a node on the Dardel supercomputer in PDC at the KTH Royal Institute of Technology, for memory reasons.
\section{Results}
In the following section, we present some experimental evaluation of our prototype method in terms of convergence of the conjugate gradient solver and interior point solver itself, as well as the conditioning of the KKT-systems and the computational performance.
\subsection{Krylov Solver Convergence}
\begin{figure}[h!]
\begin{subfigure}[t]{.45\linewidth}
    \centering
    \includegraphics[width=\linewidth]{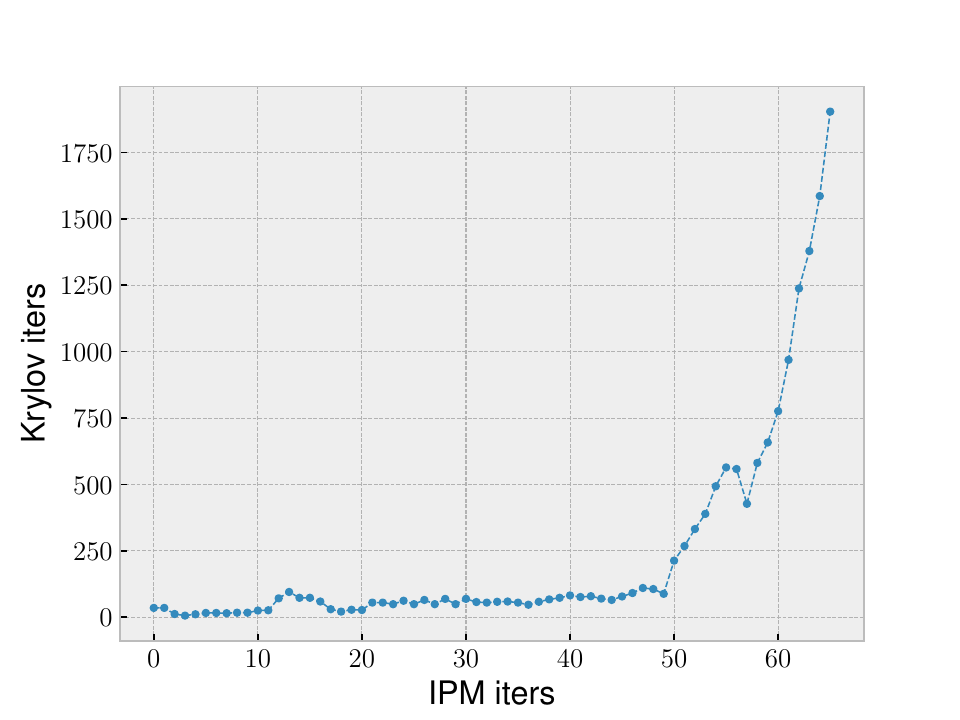}
    \caption{Photon Head and Neck case}
    \label{fig:VMAT_Krylov}
\end{subfigure}
\hfill
\begin{subfigure}[t]{.45\linewidth}
    \centering
    \includegraphics[width=\linewidth]{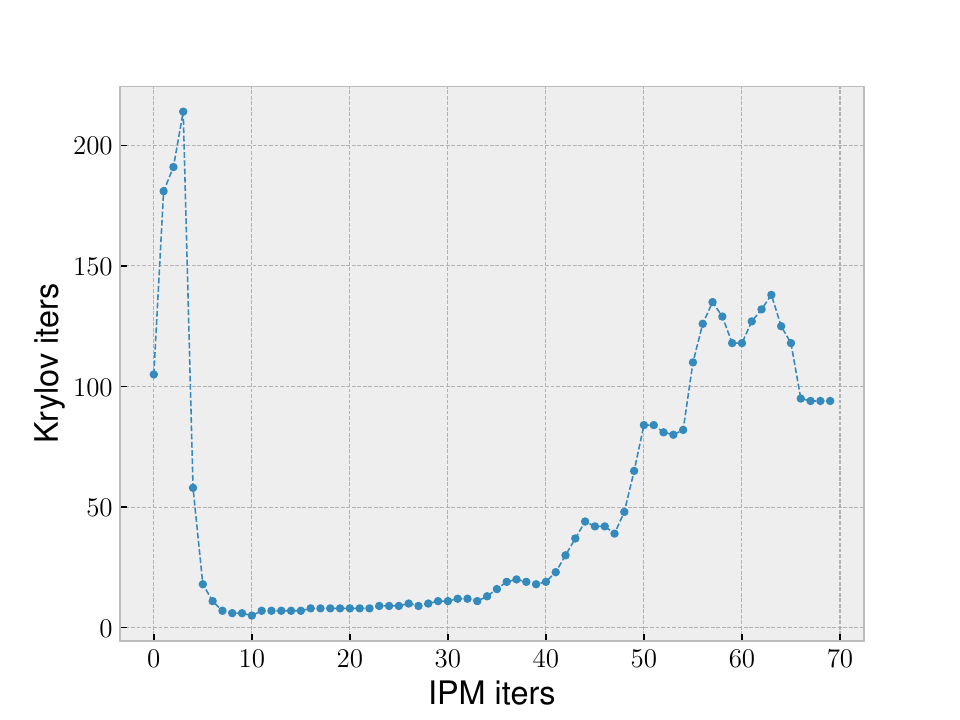}
    \caption{Proton Head and Neck case}
    \label{fig:ProtonArc_Krylov}
\end{subfigure}

\begin{subfigure}[t]{.45\linewidth}
    \centering
    \includegraphics[width=\linewidth]{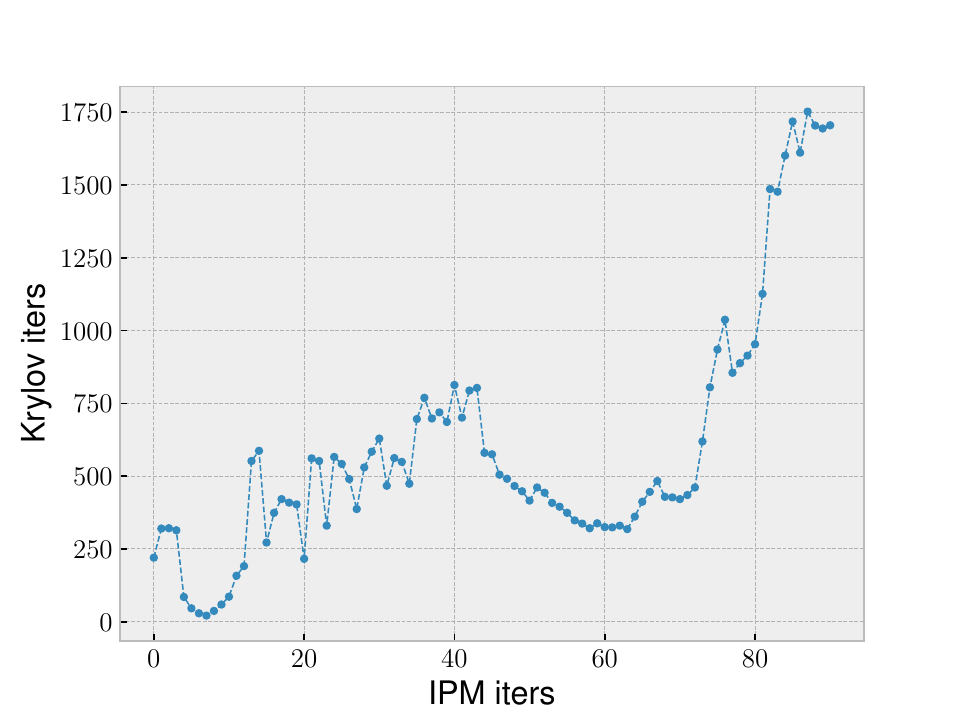}
    \caption{Proton Liver Case}
    \label{fig:ProtonLiver_Krylov}
\end{subfigure}
\hfill
\begin{subfigure}[t]{.45\linewidth}
    \centering
    \includegraphics[width=\linewidth]{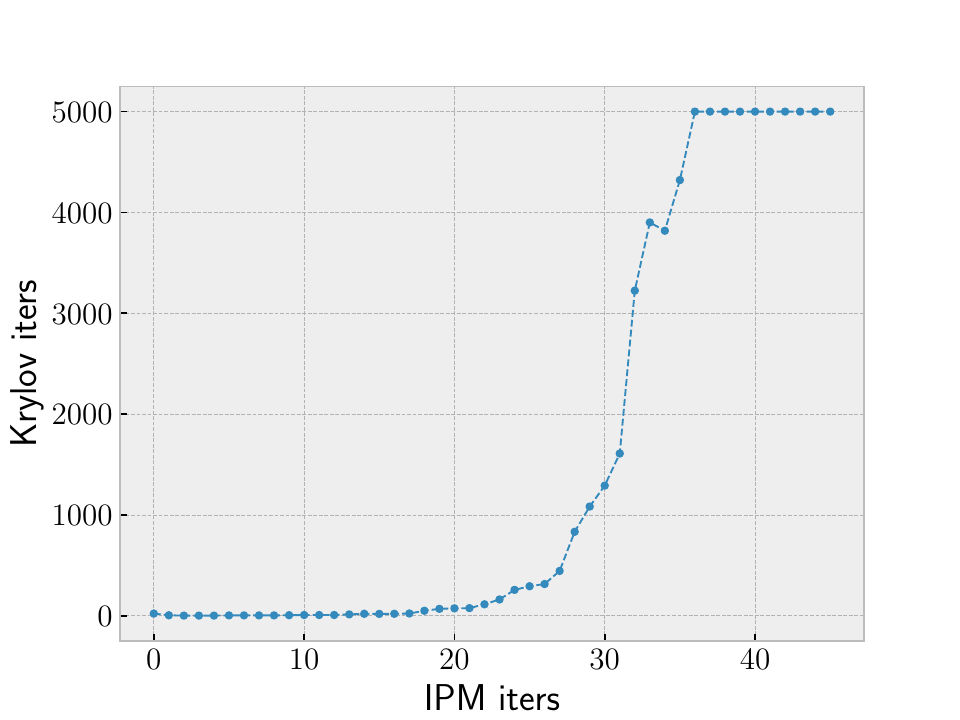}
    \caption{SVM a1a}
\end{subfigure}
\caption{Number of CG iterations at each IPM iteration for three different test problems.}
\label{fig:krylov_iters}
\end{figure}
\noindent
Fig.~\ref{fig:krylov_iters} shows the number of CG iterations required for the linear solver to converge within each IPM iteration for our test problems. Note that the maximum number of CG iterations was set to 5000 for each of the test cases, with a convergence tolerance of $10^{-7}$ for the unpreconditioned residual. As a trend, we see that all problems show a sharp increase in the number of CG iterations towards the later IPM iterations, which is consistent with the observation that the ill-conditioning of the systems arises when the barrier parameter $\mu$ gets close to zero. The proton head and neck case stands out in that it has a a spike in CG at the beginning of the optimization as well, which is also seen in the estimated condition numbers being large in the beginning in Fig.~\ref{fig:condest}. To note is that the proton head and neck case is the only one considered without general linear constraints (it has only variable bound constraints).

\subsection{IPM Solver Convergence}
\begin{figure}[h!]
\begin{subfigure}[t]{.48\linewidth}
    \centering
    \includegraphics[width=\linewidth]{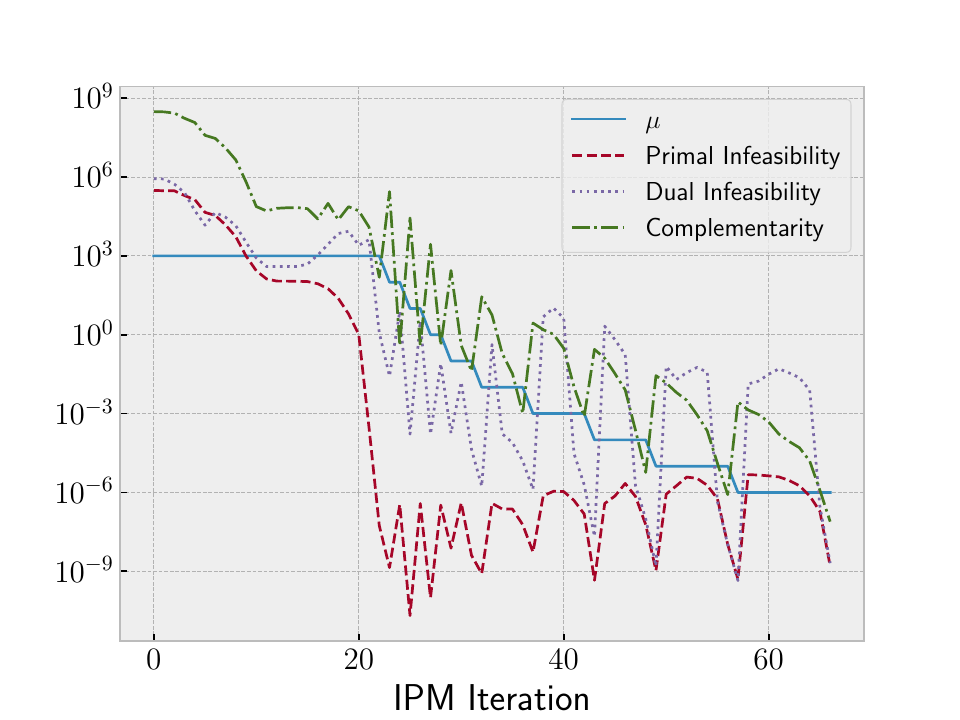}
    \caption{Photon Head and Neck case}
    \label{fig:VMAT_infeas}
\end{subfigure}
\hfill
\begin{subfigure}[t]{.48\linewidth}
    \centering
    \includegraphics[width=\linewidth]{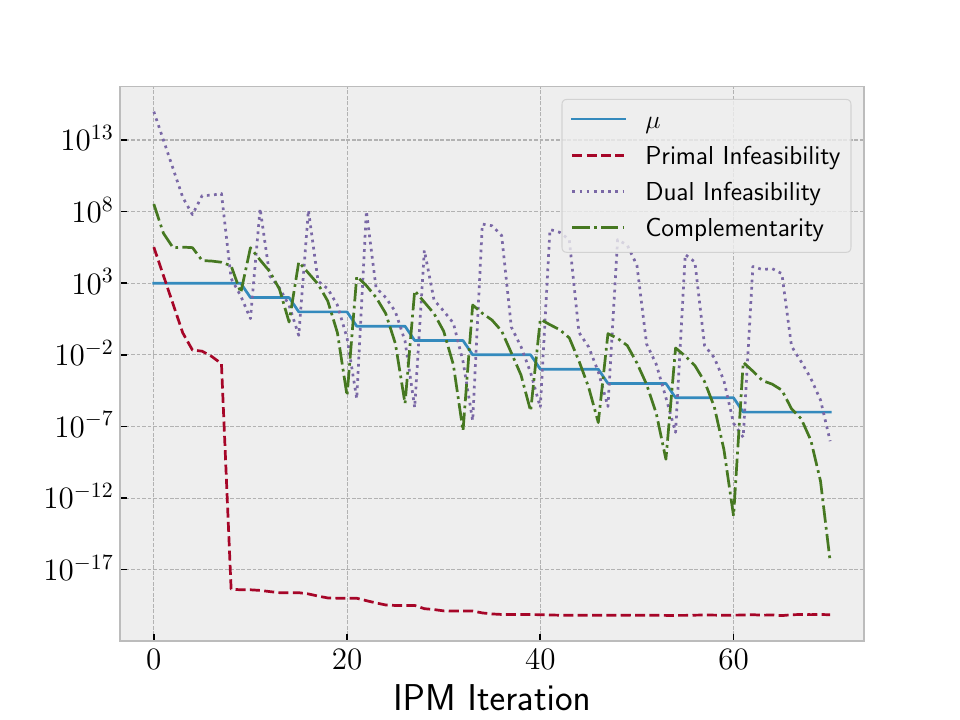}
    \caption{Proton Head and Neck case}
    \label{fig:ProtonArc_infeas}
\end{subfigure}

\begin{subfigure}[t]{.48\linewidth}
    \centering
    \includegraphics[width=\linewidth]{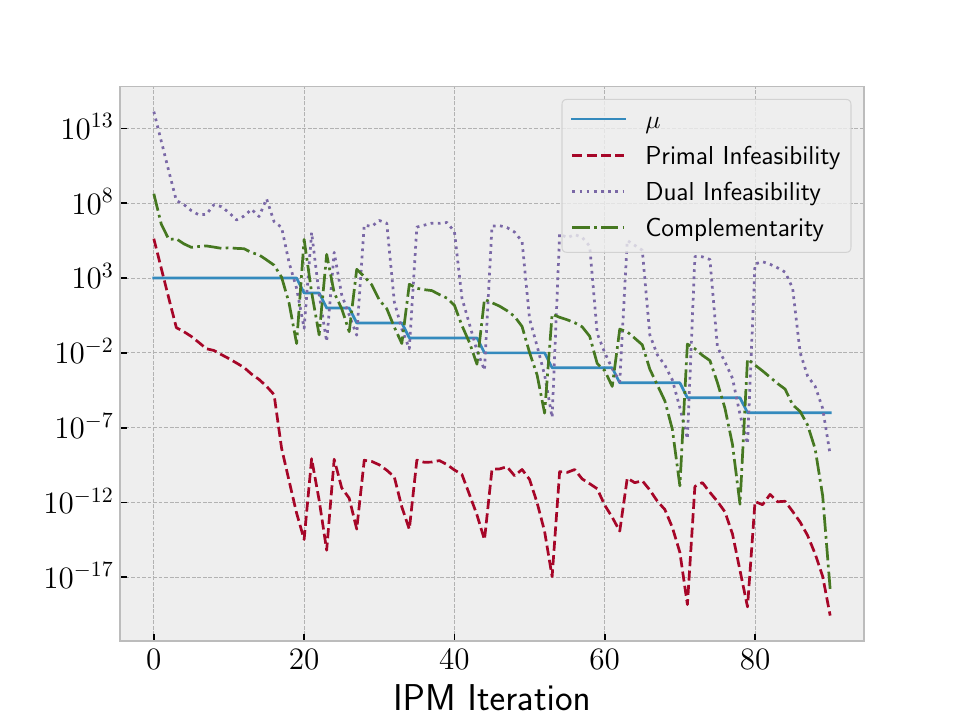}
    \caption{Proton Liver Case}
    \label{fig:ProtonLiver_infeas}
\end{subfigure}
\hfill
\begin{subfigure}[t]{.48\linewidth}
    \centering
    \includegraphics[width=\linewidth]{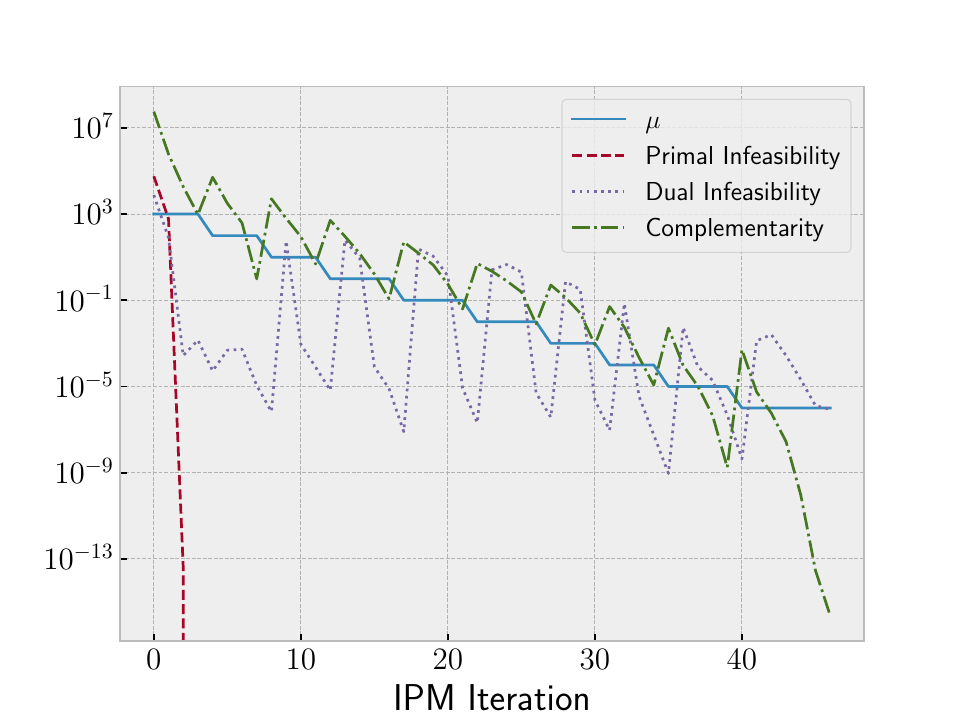}
    \caption{SVM a1a}
    \label{fig:svm_infeas}
\end{subfigure}
\caption{Convergence of our solver in terms of primal, dual and complementarity infeasibility. The value of the barrier parameter $\mu$ is shown as well, which is successively decreased as optimization progresses.}
\label{fig:infeas}
\end{figure}
\noindent
Another interesting aspect to consider is the convergence of the interior point method as a whole. In optimization solvers, the convergence is often measured with respect to the \emph{primal}, \emph{dual} and \emph{complementarity} infeasibility (among others).  The primal infeasibility is the (Euclidean) norm of
$\begin{pmatrix}
        r_{l_A} & r_{u_A} & r_{l_x} & r_{u_x}
\end{pmatrix}^T$,
the dual infeasibility is the norm of $r_H$ and the complementarity infeasibility the norm of $\begin{pmatrix} r_{c_1} & r_{c_2} & r_{c_3} & r_{c_4} \end{pmatrix}^T$.
In other words, the primal infeasibility measures the error with respect to satisfying the constraints, the dual infeasibility measures the error in stationarity of the Lagrangian, and the complementarity infeasibility the error with respect to the perturbed complementary slackness condition.

In Fig~\ref{fig:infeas}, we show the convergence in terms of the primal, dual and complementarity infeasibility over IPM iterations for our test-problems. From the figures, we see that the convergence towards optimality is far from monotonous, with the spikes in the infeasibility norms coinciding with the points when the barrier parameter $\mu$ is decreased in the solver. This could indicate that the update of $\mu$ is too aggressive. The reason could be that we use a relatively crude update rule for $\mu$ in our prototype implementation, and more sophisticated methods may give better performance. Another interesting observation from the infeasibility plots is that the dual infeasibility exhibits slower convergence compared to the complementarity and primal infeasiblities, especially for the proton cases.

Speculatively, one can observe that the dual infeasibility $r_H$, appears only once in the right-hand side of \eqref{eq:2x2_augmented}. The block in the RHS in which the dual infeasibility appears is also contains multiple terms scaled by the inverse of the slack variables, which may be large for \emph{active} constraints. The doubly augmented system \eqref{eq:doubly_augmented} introduces an additional term in the top block of the RHS. In view of Krylov methods as algorithms seeking least-norm solutions in a given Krylov subspace for a set of linear system of equations, this may give a partial explanation why the dual infeasibility lags behind in our case, as the contribution to the RHS in the linear system from the corresponding term can be relatively small.
\subsection{Numerical Stability and Conditioning of KKT System}
\begin{figure}[h!]
    \centering
    \includegraphics[width=.6\linewidth]{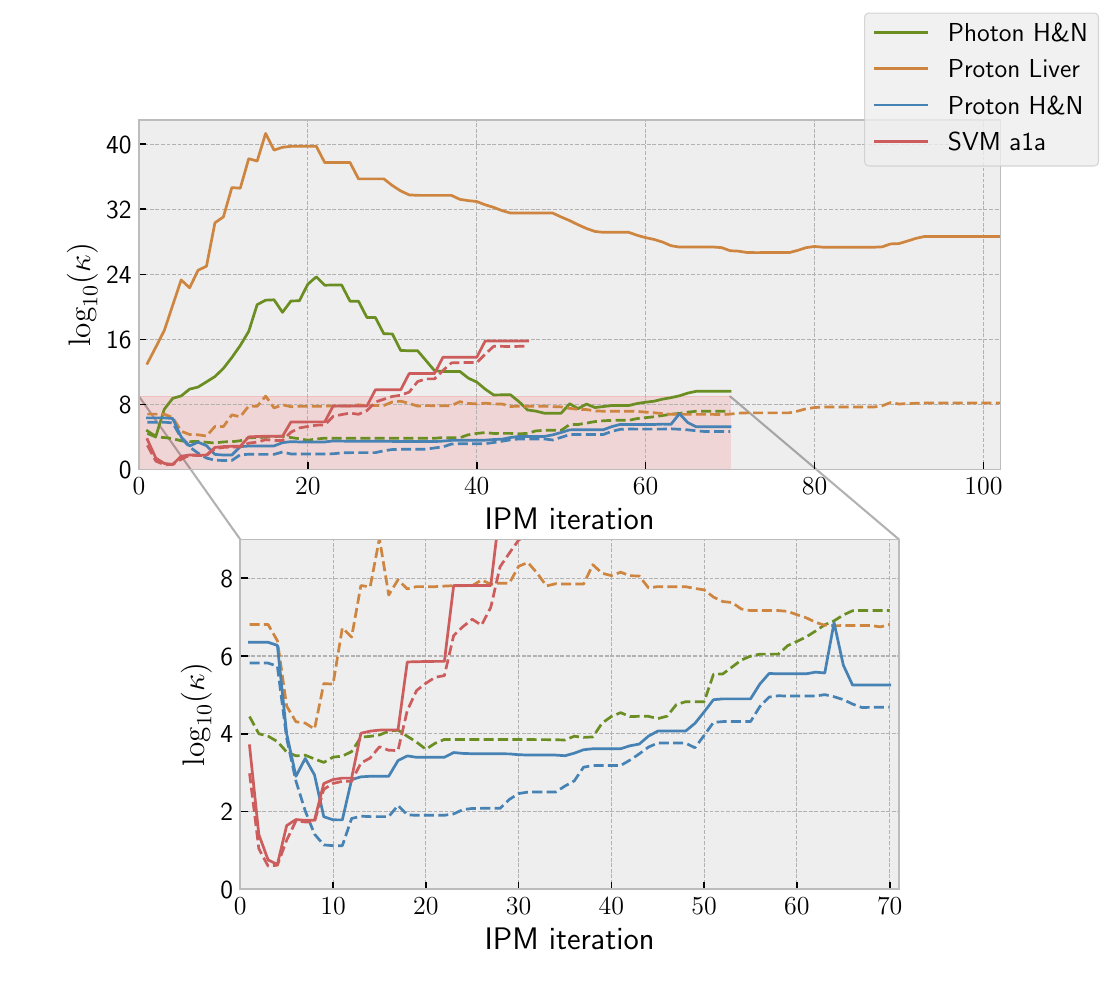}
    \caption{Condition number estimates using Matlab's \texttt{condest} throughout the optimization iterations. The bottom plot is zoomed in on the red region (with the Photon H\&N line removed for clarity). The dashed lines are condition numbers after Jacobi preconditioning.}
    \label{fig:condest}
\end{figure}
\noindent
As discussed previously, one of the main challenges in using iterative linear solvers for interior point methods is the structured ill-conditioning in the linear systems. To study how this conditioning affects our problem, we evaluate how the condition number $\kappa$ of the doubly augmented KKT system \eqref{eq:doubly_augmented} that we solve in each iteration changes throughout the optimization. The condition numbers of the resulting matrices are estimated using Matlab's \texttt{condest} function, which we modify slightly by using Cholesky instead of LU-factorization internally, since our matrices are symmetric positive definite. \texttt{condest} gives an estimate of the condition number in the $L1$-norm and is based on an algorithm proposed by Higham in \cite{higham1988fortran}.

Fig.~\ref{fig:condest} shows the results of the condition number analysis. The solid lines show the un-preconditioned condition numbers, with the dashed lines showing the condition numbers with Jacobi preconditioning. The un-preconditioned KKT systems for the Photon H\&N and Proton Liver cases show extreme ill-conditioning, especially in the middle of the optimization, with estimated condition numbers up to the order $10^{41}$ for the Proton Liver case and $10^{23}$ for the Photon H\&N case. While the accuracy of Matlab's estimation using \texttt{condest} at such extreme ill-conditioning may be questioned, suffice it to say that the un-preconditioned matrices are close to singular. However, we see that the Jacobi preconditioning does manage to improve the conditioning of those matrices significantly, reducing the condition number to around $10^8$ (or less) for both cases.


\subsection{Performance Analysis}
\begin{figure}[h!]
    \centering
    \includegraphics[width=.5\linewidth]{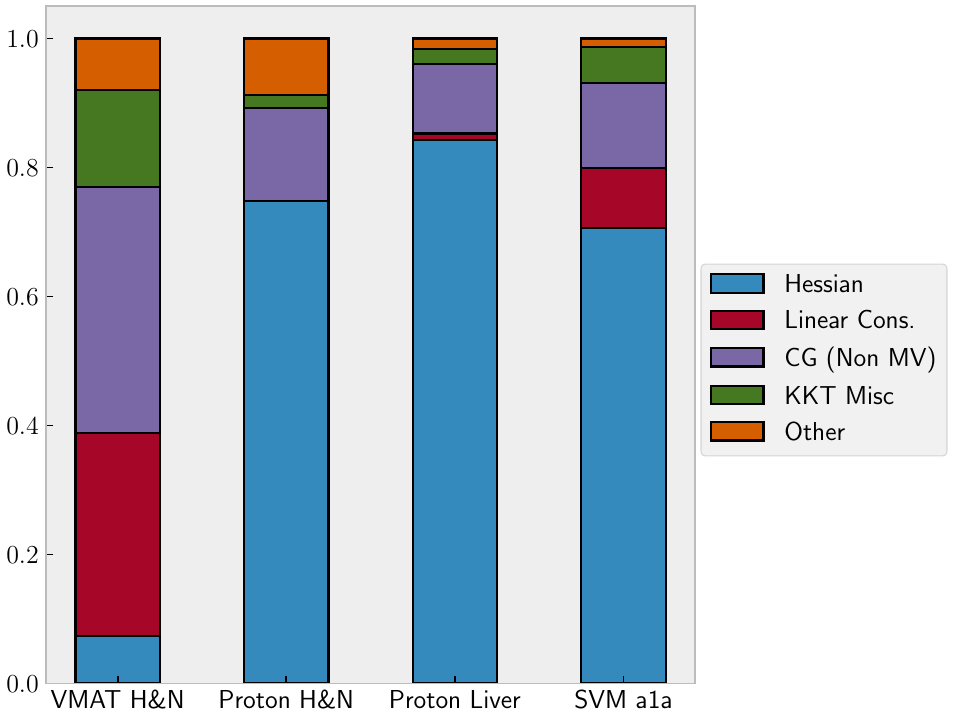}
    \caption{Relative (normalized) run-time spent in different computational kernels for our solver.}
    \label{fig:profiling}
\end{figure}
\begin{table}[h!]
    \centering
    \caption{Total run times in seconds spent in different parts of the solver for the different problems. CG (Not MV) is the time spent in the conjugate gradients solver excluding matrix-vector products. The relative times in different parts are also visualized in Fig.~\ref{fig:profiling}.}
    \label{tab:profiling}
    \begin{tabular}{|c|c|c|c|c|c|c|}
    \hline
        Case & Total & Hessian & Lin. Cons. & CG (Not MV) & KKT Misc. & Other \\
    \hline
        Photon H\&N & 13.6 & 0.997 & 4.30 & 5.20 & 2.06 & 1.08 \\
        Proton H\&N & 6.28 & 4.70 & 0.002 & 0.907 & 0.127 & 0.551 \\
        Proton Liver & 324 & 222 & 62.8 & 28.4 & 6.13 & 4.29 \\
        SVM a1a & 2.78 & 1.96 & 0.256 & 0.367 & 0.157 & 0.035 \\
    \hline
    \end{tabular}
\end{table}
\noindent
Fig.~\ref{fig:profiling} shows the relative time spent in different parts of the code for the three test problems considered in this work. The absolute run-times for the different parts is given in Table~\ref{tab:profiling}. The parts we included in the profiling are: Matrix-vector multiplication with the quasi-Newton Hessian, matrix-vector multiplication with the constraint matrix, remaining time in the CG solver (excluding matrix vector multiplication) remaining time spent in the matrix-vector multiplication with the doubly augmented KKT system in \eqref{eq:doubly_augmented}, and finally the "Other" category comprising the remaining time spent in other parts of the solver.

From Fig.~\ref{fig:profiling}, we see that for both proton cases and the SVM problem, the solver spends the majority of the time in computing (dense) matrix-vector products with the Hessian, while for the photon radiation therapy problem, a significant amount of time is spent in the CG solver itself, as well as for multiplication with the linear constraints. Overall, we believe the performance analysis shows that there is potential for improved performance when moving to GPU, especially for the proton radiation therapy problems and SVM problems.



\section{Related Work}
The topic of iterative linear solvers in interior point methods has attracted much research, which we briefly summarize in the following section. Preconditioners for KKT systems in interior point method have been studied extensively previously, for instance in \cite{bergamaschi2004preconditioning,karim2022efficient,gill1992preconditioners,rees2007preconditioner,cui2019implementation,zilli2022block} among many others. A general overview of HPC in the space of optimization and optimization software can be found in \cite{liu2022survey}. The topic of parallel computing in optimization and operations research in general has also been surveyed previously \cite{schryen2020parallel}.

Practical studies where iterative linear solvers are used for different kinds of optimization problems can be found in \cite{regev2022hykkt}, where a type of hybrid direct-iterative solution method is evaluated on very large problems in optimal power flow. In the context of interior point methods for linear programming (LP), preconditioned Krylov methods were studied in e.g.\cite{cui2019implementation,chowdhury2022faster}.
\section{Conclusions and Future Work}
In this work, we have presented our prototype interior point method for quadratic programming that uses an iterative linear solver for the KKT systems arising in each iteration. We demonstrate that the method can solve real optimization problems from radiation therapy to acceptable levels of accuracy and within reasonable time. From analyzing the performance of our implementation using tracing and profiling, we believe that our method is suitable for GPU acceleration, which we will investigate further in future work. Overall, we believe that interior point methods using Krylov solvers give a promising path forward for GPU accelerated interior point methods, which hold great promise for e.g. computational efficieny in treatment plan optimization for radiation therapy.

While we believe this initial study on the use of iterative linear solvers for optimization of radiation therapy treatment plans shows promise, there are many interesting questions and problems remaining for future research. Among those is the porting of the code to be able to run on GPU accelerators, looking at improved preconditioners for the KKT systems and investigating the applicability of the method to optimization problems from other domains.
\section{Acknowledgment}
The computations were enabled by resources provided by the National Academic Infrastructure for Supercomputing in Sweden (NAISS) at PDC, partially funded by the Swedish Research Council through grant agreement no. 2022-06725.
\bibliographystyle{splncs04}
\bibliography{References}
\end{document}